\documentclass[11pt,fleqn]{article}

\usepackage{latexsym}

\usepackage{amsmath}
\usepackage{amsthm}
\usepackage{amssymb}
\usepackage{amsfonts}

\newcommand{\be}{\begin{equation} } 
\newcommand{\ee}{\end{equation} \par \noindent}

\newcommand{\rf}[1]{(\ref{#1})}

\newcommand{\ba}{\begin{array}}
\newcommand{\ea}{\end{array}}
\newcommand{\const}{{\rm const}}

\newcommand{\R}{{\mathbb R}}
\newcommand{\Z}{{\mathbb Z}}

\newcommand{\T}{{\mathbb T}}

\newcommand{\rr}{{\vec r}}
\newcommand{\nn}{{\vec n}}

\newtheorem{prop}{Proposition}

\newtheorem{rem}{Remark}
\newtheorem{cor}{Corollary}
\newtheorem{Def}{Definition}

\newenvironment{Proof}{\par \vspace{2ex} \par
\noindent \small {\it Proof:}}{\hfill $\Box$ 
\vspace{2ex} \par }

\begin{document}

\title{\bf 
Pseudospherical surfaces on time scales}

\author{{\bf
Jan L.\ Cie\'sli\'nski\thanks{
E-mail: \tt janek\,@\,alpha.uwb.edu.pl}}
%\\ {\bf Iwona \'Swis\l ocka}
\\ {\footnotesize Uniwersytet w Bia\l ymstoku, Instytut Fizyki Teoretycznej}
\\ {\footnotesize ul.\ Lipowa 41, 15-424  Bia\l ystok, Poland} 
}

\maketitle

\begin{abstract} 
We define and discuss the notion of pseudospherical  surfaces in asymptotic coordinates on time scales. Two special cases, namely dicrete pseudospherical surfaces and smooth pseudosperical surfaces are consistent with this description. In particular, we define the Gaussian curvature in the discrete case. 
\end{abstract}

\par \vspace{0.5cm} \par
\par \vspace{0.5cm} \par
\par \vspace{0.5cm} \par

\section{Introduction}

A time scale is an arbitrary nonempty closed subset of the real numbers \cite{Hi}. Typical examples are $\R$ and $\Z$. The time scales were introduced in order to unify differential and difference calculus \cite{Hi,Hi2}. 
Partial differentiation, tangent lines and tangent planes on time scales have been introduced recently \cite{BG}.

On the other hand, besides the differential geometry, there exists also the difference geometry \cite{Sau}. In the last years one can observe a fast development of the integrable difference geometry (see, for instance, \cite{BS,Ci-fam,CDS,Do,DS-nets,WKS}) related closely to the classical differential geometry based on explicit constructions and transformations \cite{DSM,Eis1}.  
A natural idea is to unify the difference and differential geometries and to formulate the integrable geometry on time scales.

In this paper we propose such formulation for pseudospherical immersions (surfaces of constant negative Gaussian curvature). The discrete pseudospherical surfaces have been introduced a long time ago \cite{Sau,Wun}, and studied intensively in the last years \cite{BP-pseudo}. 
The discrete pseudospherical surfaces have been recently expressed in terms of time scales \cite{Sw-mgr}. However, the assumption that all points are isolated was essential and the Gaussian curvature was not discussed at all. In the present paper we use a different starting point  and all cases are described in a unified framework.

\section{Differentiation on time scales}

This section collects basic notions and results concerning the differential calculus on time scales, compare \cite{BG}. To avoid some unimportant complications we confine ourselves to time scales which are not bounded neither from above nor from below
(more general case is considered in \cite{BG}, as far as partial derivatives are concerned).

\begin{Def}[\cite{Hi}] Let a time scale $\T$ is given. The maps $\sigma : \T \rightarrow \T$ and $\rho : \T \rightarrow \T$, defined by 
\be  \ba{l}
\sigma (u) := \inf \{ v \in \T : v > u \}  \ ,  \\[2ex]
\rho (u) := \sup \{ v \in \T : v < u \}  \ ,
\ea \ee
are called jump operator and backward jump operator, respectively.
\end{Def}

\begin{Def}[\cite{Hi}] The points  $u \in \T$ can be classified as follows 
\begin{itemize} 
\item \quad $\sigma (u) > u$ \qquad $\Longrightarrow$ \quad $u$ is right-scattered , 
\item  \quad $\sigma (u) = u$  \qquad $\Longrightarrow$ \quad $u$ is right-dense ,
\item \quad $\rho (u) < u$  \qquad $\Longrightarrow$ \quad $u$ is left-scattered ,
\item \quad $\rho (u) = u$  \qquad $\Longrightarrow$ \quad $u$ is left-dense ,
\item \quad $\rho (u) < u < \sigma (u)$  \quad  $\Longrightarrow$ \quad $u$ is isolated . 
\end{itemize}
\end{Def}

\begin{Def}[\cite{BG}]
The delta derivative is defined as
\be
\frac{\partial f (t) }{\Delta t} = \lim_{\stackrel{s \rightarrow t}{s \neq \sigma (t)}} \frac{ f ( \sigma (t) ) - f ( t )}{\sigma (t) - s} \ ,
\ee
and the nabla derivative is defined by
\be
\frac{\partial f (t) }{\nabla t} = \lim_{\stackrel{s \rightarrow t}{s \neq \rho (t)}} \frac{ f ( \rho (t) ) - f ( t )}{\rho (t) - s} \ .
\ee
\end{Def}

\begin{Def}[\cite{BG}]  \label{complete}
 We say that a function $f : \T \rightarrow \R$ is completely delta differentiable at a point $t_0 \in \T$, if there exist a number $A$ such that 
\[  \ba{l}
f (t) - f(t_0) = A (t - t_0) + (t - t_0) \ \alpha (t_0, t) \ , \\[2ex]
f(t) - f(\sigma (t_0)) = A ( t - \sigma (t_0) ) + (t - \sigma (t_0)) \ \beta (t_0, t)
\ ,
\ea  \]
where \ $\alpha (t_0, t_0) = 0$, $\beta (t_0, t_0) = 0$, $\displaystyle \lim_{t\rightarrow t_0} \alpha (t_0, t) = 0$, and \ $\displaystyle  \lim_{t\rightarrow t_0} \beta (t_0, t) = 0$. 

\end{Def}

\begin{prop}[\cite{BG}]
If the function $f$ is completely delta differentiable at $t_0$, then the graph of this function has the uniquely determined delta tangent line at the point $P_0 = (t_0, f(t_0))$ specified by the equation
\[
y - f(t_0) = \frac{\partial f (t_0) }{\Delta t} (x - t_0) 
\] 
\end{prop}

In this paper we fix our attention on functions defined on two-dimensional time scales, i.e., on $\T_1 \times \T_2 $, where $\T_1, \T_2$ are given time scales. 
The extension on $n$-dimensional time scales is usually straightforward. We denote: 
\be  \ba{l}
t \equiv (t_1, t_2) \in \T_1 \times \T_2 \ ,  \\[2ex]
\sigma_1 (t) = ( \sigma (t_1) , t_2)  \ , \quad \sigma_2 (t) = (t_1, \sigma (t_2)) \ , \\[2ex]
\rho_1 (t) = ( \rho (t_1) , t_2)  \ , \quad \rho_2 (t) = (t_1, \rho (t_2)) \ .
\ea \ee

\begin{rem}
In the discrete case ($\T_1 = \T_2 = \Z$)  we have $\sigma_j (u) = T_j u$ and $\rho_j (u) = T_j^{-1} u$, where $T_j$ mean usual shift operators. Therefore delta and nabla differentiation can be associated with forward and backward data, respectively \cite{Do}. 
\end{rem}

\begin{Def}[\cite{BG}]
The partial delta derivative is defined as
\be
\frac{\partial f (t) }{\Delta t_j} = \lim_{\stackrel{s_j \rightarrow t_j}{s_j \neq \sigma (t_j)}} \frac{ f ( \sigma_j (t) ) - f ( t )}{\sigma (t_j) - s_j} \ .
\ee
The definition of the partial nabla derivative is analogical.
\end{Def}

In the continuous case (e.g., $\T_1 = \T_2 = \R$) the delta derivative coincides with the right-hand derivative, while the nabla derivative coincides with the left-hand derivative. Note that all results and definitions in terms of delta derivatives have their nabla derivatives analogues. 

\begin{prop}[\cite{BG}]
If the mixed partial delta derivatives exist in a neighbourhood of $t_0 \in \T_1 \times \T_2$  and are continuous at $t = t_0$, then 
\[
   \frac{\partial^2 f (t_0) }{\Delta t_1  \Delta t_2 } = \frac{\partial^2 f (t_0) }{\Delta t_2 \Delta t_1 } \ .
\]
\end{prop}

The definition of the complete delta differentiability is similar to Definition~\ref{complete}, see \cite{BG}, Definition 2.1. Instead of this definition we present here an important sufficient condition. 

\begin{prop}[\cite{BG}]
Let a function $f : \T_1 \times \T_2 \rightarrow \R$ be continuous and have first order partial derivatives in a neighbourhood of $t_0$. If these derivatives are continuous at $t_0$, then $f$ is completely delta differentiable at $t_0$.
\end{prop}

\begin{Def}[\cite{BG}] Let $z = f (x, y)$ ($x \in \T_1$, $y \in \T_2$) be a given surface (on the time scale) in $\R^3$.  
A plane $\Omega_0$ passing through  $P_0 = (t_0, s_0, f((t_0, s_0))$ (where $t_0 \in \T_1$, $s_0 \in \T_2$) is called the delta tangent plane to the surface $S$ at the point $P_0$ if 
\begin{enumerate}
\item  $\Omega_0$ passes also through the points $P_0^{\sigma_1} = (\sigma_1 (t_0), s_0, f(\sigma_1(t_0), s_0))$ and $P_0^{\sigma_2} = (t_0, \sigma_2 (s_0), f(t_0, \sigma_2 (s_0))$; 
\item  if $P_0$ is not isolated point of $S$ then 
\[ 
\lim_{\stackrel{P\rightarrow P_0}{P\neq P_0}} \frac{d(P,\Omega_0)}{d(P,P_0)} = 0 \ ,
\]
where $P$ is a moving point of the surface $S$, $d(P,\Omega_0)$ is the distance from $P$ to the plane $\Omega_0$ and $d(P,P_0)$ is the distance between $P$ and $P_0$.
\end{enumerate}

\end{Def}

Delta tangent line is defined in an analogous way. If $P_0$ is an isolated point of the curve $\Gamma$ (hence $P_0 \neq P_0^\sigma$), then 
the delta tangent line to $\Gamma$ at $P_0$ coincides with the unique line through the points $P_0$ and $P_0^\sigma$. 

Similarly, if $P_0 \neq P_0^{\sigma_1}$ and $P_0^{\sigma_2} \neq P_0$ (hence also $P_0^{\sigma_1} \neq P_0^{\sigma_2}$), then the delta tangent plane to the surface $S$ at $P_0$ (if exists) coincides with the unique plane through $P_0, P_0^{\sigma_1}$ and $P_0^{\sigma_2}$. 

\begin{prop}[\cite{BG}]   \label{tanplane}
If the function $f : \T_1 \times \T_2 \rightarrow {\mathbb R}$ is completely delta differentiable at $(t_0,s_0)$, then the surface represented by this function has the uniquely determined delta tangent plane at the point $P_0 = (t_0, s_0, f(t_0, s_0))$ specified by the equation 
\be
z = f(t_0,s_0) + \frac{\partial f (t_0, s_0)}{\Delta t} (x - t_0) + \frac{\partial f (t_0, s_0)}{\Delta s} (y - s_0) 
\ee
where $(x,y,z)$ is the current point of the plane.
\end{prop}

In the following sections of this paper we define pseudospherical surfaces on time scales in terms of delta derivatives. In order to simplify the notation the delta derivatives will be denoted by 
\be  \label{D}
D_j f \equiv \frac{\partial f (t) }{\Delta t_j} \ .
\ee
Proposition~\ref{tanplane} suggests that in  geometrical contexts it is more natural to use complete delta differentiability rather than delta differentiability.

\section{Pseudospherical surfaces}
Let us consider a surface immersed in $\R^3$ explicitly described by a position vector $\rr = \rr (u,v)$. Denoting the normal vector by $\nn$ we define the so called fundamental forms:
\[
  I := d \rr \cdot d \rr \ , \qquad II := - d \rr \cdot d \nn \ ,
\]
where the center dot denotes the standard scalar product in $\R^3$. We denote the coefficients of the fundamental forms in a traditional way:
\be  \label{fund}  \ba{l}
 I = E du^2 + 2 F du dv + G dv^2 \ , \\[2ex]
II = L du^2 + 2 M du dv + N dv^2 \ .
\ea \ee 
Hence,
\be  \ba{c} \label{contform}
E = \rr_u \cdot \rr_u \ , \qquad F = \rr_u \cdot \rr_v \ , \qquad G = \rr_v \cdot \rr_v \ , \\[2ex]
L = - \nn_u \cdot \rr_u \ , \qquad M = - \nn_u \cdot \rr_v \ , \qquad N = - \nn_v \cdot \rr_v \ .
\ea \ee
The Gaussian curvature $K$ and the mean curvature $H$ are given by:
\be \label{KH} 
K = \frac{L N - M^2}{W^2} \ , \qquad
H = \frac{E N - 2 F M + G L}{2 W^2} \ ,
 \ee
where $W = E G - F^2$ (by assumption $W \neq 0$, i.e., the first fundamental form is not degenerated). 
The coefficents $E, F, G, L, M, N$ satisfy the Gauss equation \cite{Ampol}
\be  \ba{c} \label{Gauss}  \displaystyle
K = - \frac{1}{4 W^2} \left| \begin{array}{ccc} E & E_u & E_v \\ F & F_u & F_v \\ G & G_u  & G_v  \end{array}  \right| + \frac{1}{2 W} \left( \frac{\partial}{\partial u} \frac{F_v - G_u}{W}  + \frac{\partial}{\partial v} \frac{F_u - E_v}{W}  \right) \ ,
\ea \ee
and two Peterson-Mainardi-Codazzi equations
\be \label{MC}   \ba{ll}  \displaystyle 
L_v  - M_u - H (E_v - F_u) + \frac{1}{2 W^2} \left| \begin{array}{ccc} E & E_u & L \\ F & F_u & M \\ G & G_u  & N  \end{array}  \right| = 0 \ ,
\\[5ex]    \displaystyle 
M_v  - N_u - H (F_v - G_u) + \frac{1}{2 W^2} \left| \begin{array}{ccc} E & E_v & L \\ F & F_v & M \\ G & G_v  & N  \end{array}  \right| = 0 \ .
\ea \ee
The Bonnet theorem says that any solution of the system \rf{Gauss},\rf{MC} 
implicitly defines a surface immersed in $\R^3$ (provided that $E > 0$ and $W > 0$) \cite{Ampol}.

\begin{prop}  \label{czeb}
Let asymptotic lines on a surface admit parameterization by Chebyshev coordinates, i.e., the fundamental forms are expressed in terms of two real functions $F, M : \R^2 \supset \Omega \rightarrow \R$ as follows  
\be  \label{fundpseudo}
 I = du^2 + 2 F(u,v) du dv + dv^2 \ , \qquad
II = 2 M (u,v) du dv \ , 
\ee 
then the surface ${\rr} = {\rr} (u,v)$ ($u,v \in \Omega)$, implicitly defined by the fundamental forms \rf{fundpseudo}  has a constant negative Gaussian curvature. 
\end{prop}

\begin{Proof} Substituting $E = G = 1$ and $L = N = 0$ to \rf{KH} and \rf{MC}, we get
\[  \ba{l}  \displaystyle
K = - \frac{M^2}{1 - F^2} \ , \quad H = - \frac{FM}{1 - F^2} \ , \quad
 M_u - H F_u = 0 \ , \quad  M_v - H F_v = 0 \ .
\ea \]
Hence,
\[
M M_u (1- F^2) + F F_u M^2 = 0 \ , \quad  M M_v (1- F^2) + F F_v M^2 = 0 \ , 
\]
which means \ $\displaystyle  \frac{M^2}{1 - F^2} = \const > 0$.
Therefore, $K = \const < 0$. 
\end{Proof}

\begin{rem}
The assumptions of the Lemma~\ref{czeb} can be rewritten as
\be
     (\rr_u )^2 = (\rr_v )^2 = 1 \ ,  \qquad  \nn_u \cdot \rr_u = \nn_v \cdot \rr_v = 0 \ , 
\ee
and the conclusion of Lemma~\ref{czeb} states
\be  \label{Kcont}
   K \equiv \frac{ - M^2}{1 - F^2} = \const < 0 \ .
\ee
\end{rem}

We recall that asymptotic lines are characterized by $L = N = 0 $, i.e., the second fundamental form is given by \rf{fundpseudo}. 
Having Chebyshev coordinates $u, v$ we can consider more general parameterization of asymptotic lines, namely: $\tilde u = f (u)$, $\tilde v = g (v)$. They are called weak Chebyshev coordinates.

\section{Discrete pseudospherical surfaces}

In the discrete case the time scale $\T_1 \times \T_2$ contains only  isolated points. We confine ourselves to the case $\T_1 = \T_2 = a \Z$, where $a$ is a fixed constant (the mesh size). 

\begin{rem} Let $\T_1 = \T_2 = a \Z$ and $f : \T_1 \times \T_2 \rightarrow \R$, then we denote
\be  \label{Delta}
   \Delta_j f =  \frac{\partial f (t)}{\Delta_j t} = \frac{T_j f - f}{a} \ .
\ee
Therefore, in the discrete case $D_j = \Delta_j$. In particular, for $a=1$ we have $\Delta_j = T_j - 1$.
\end{rem}

The discrete analogue of pseudospherical surfaces endowed with Chebyshev coordinates \rf{fundpseudo}, i.e., discrete Chebyshev net, is defined as follows,  compare \cite{Wun}.

\begin{Def}[\cite{BP-pseudo}] \label{pseudodis} Discrete Chebyshev net (discrete $K$-surface) is an  immersion $\rr: a\Z \times a \Z \ni (a m, a n) \rightarrow \rr (a m, a n) \in \R^3$ such that for any $m,n$
\begin{itemize}
\item $|\Delta_1 \rr| = |\Delta_2 \rr| = 1$ \ ,
\item the points $\rr$, $T_1 \rr$, $T_2 \rr$, $T_1^{-1} \rr$, $T_2^{-1} \rr$  are  coplanar (we denote this plane by $\pi (\rr)$). 
\end{itemize}
\end{Def}

\noindent By the discrete immersion we mean that $\Delta_1 \rr$ and $\Delta_2 \rr$ are linearly independent for any $m,n$.

 Similarly one can discretize weak Chebyshev coordinates \cite{BP-pseudo}. However, in this paper we confine ourselves only to discrete Chebyshev nets.

The plane $\pi (\rr)$ can be interpreted, obviously, as the discrete analogue of the tangent plane. Therefore  
\be
  \nn := \frac{ \Delta_1 \rr \times \Delta_2 \rr }{| [\Delta_1 \rr, \Delta_2 \rr] |} = \frac{ \Delta_1 \rr \times \Delta_2 \rr }{\sqrt{ 1 - \Delta_1 \rr \cdot \Delta_2 \rr } } \ ,
\ee
is the discrete analogue of the normal vector (here the cross means the vector product). 

\begin{prop}
In the discrete case 
\[
\Delta_1 \nn \cdot \Delta_1 \rr = 0 \quad \Longleftrightarrow \quad \Delta_1 \rr \ , \ T_1 (\Delta_1 \rr)\ , \ T_1 (\Delta_2 \rr) \ {\rm are \ coplanar}  .
\]
\[
\Delta_2 \nn \cdot \Delta_2 \rr = 0 \quad \Longleftrightarrow \quad \Delta_2 \rr \ , \ T_2 (\Delta_1 \rr)\ , \ T_2 (\Delta_2 \rr) \ {\rm are \ coplanar}  .
\]
\end{prop}

\begin{Proof} From the definition of $\nn$ it follows: $\nn \cdot \Delta_1 \rr = 0$, $T_1 \nn \cdot T_1 \Delta_1 \rr = 0$ and $T_1 \nn \cdot T_1 \Delta_2 \rr = 0$. Then $\Delta_1 \nn \cdot \Delta_1 \rr = 0  \ \Longleftrightarrow \ T_1 \nn \cdot \Delta_1 \rr = \nn \cdot \Delta_1 \rr $. Hence, $T_1 \nn \cdot \Delta_1 \rr = 0$. Therefore, 
$\Delta_1 \rr$, $T_1 \Delta_1 \rr$ and $T_1 \Delta_2 \rr$ are co-planar.
The proof of the second statement is similar.
\end{Proof}

\begin{cor} \label{cordis} In the discrete case  
$\rr$,   $T_1 \rr$, $T_2 \rr$, $T_1^{-1} \rr$, $T_2^{-1} \rr$ \ are  coplanar if and only if \ 
$\Delta_1 \nn \cdot \Delta_1 \rr = 0$ and  
$\Delta_2 \nn \cdot \Delta_2 \rr = 0$. 
\end{cor}

In the next part of this section we consider the tetrahedron $ABCD$:
\[
\rr \equiv A \ , \quad  T_1 \rr \equiv B \ , \quad T_2 \rr \equiv D \ , \quad T_1 T_2 \rr \equiv C \ .
\]
The angle between $\Delta_1 \rr$ and $\Delta_2 \rr$ will be denoted by $\phi$ and the angle between $-\Delta_2 \rr$ and $T_2 \Delta_1 \rr$ will be denoted by $\psi$. The tetrahedron ABCD is uniquely defined by specifying $a, \phi, \psi$.

\begin{prop}  \label{tconst}
The angle $\theta_j$ between $\pi (\rr)$ and $T_j ( \pi (\rr))$ $(j=1,2)$ is constant, i.e., $\theta_j (m,n) = \theta = \const$. 
\end{prop}

\begin{Proof} 
The transformation \ $T_1 \rr \longleftrightarrow T_2 \rr$ \ is an isometry of the tetrahedron $ABCD$. Hence the angle between $\pi (\rr)$, $\pi (T_1 \rr)$ is equal to the angle between $\pi (\rr) $,  $\pi (T_2 \rr)$. 
The transformation \ $\rr \longleftrightarrow T_1 T_2 \rr$, \ $T_1 \rr \longleftrightarrow T_2 \rr$ \ is another isometry of this tetrahedron. Thus \ $\pi(\rr) \longleftrightarrow \pi(T_1 T_2 \rr)$, \ $\pi (T_1 \rr) \longleftrightarrow \pi (T_2 \rr)$. Hence, the angle between $\pi (\rr)$,  $\pi (T_j \rr)$ is equal to the angle between $\pi (T_k \rr)$, $\pi (T_k T_j \rr)$, which means that this angle does not depend on $m,n$. 
\end{Proof}

\begin{prop} \label{Kprop}
In the discrete case $K$ defined by
\be  \label{Kdisc}
 K = - \frac{ (\Delta_1 \nn\cdot \Delta_2 \rr ) (\Delta_2 \nn \cdot \Delta_1 \rr )   }{ 1 - (\Delta_1 \rr \cdot \Delta_2 \rr )^2 }  
\ee
is constant (i.e., does not depend on $m,n$). Moreover 
\be
  K = - \frac{\sin^2 \theta}{a^2} \ .
\ee
\end{prop}

\begin{Proof} 
Taking into account $|\vec{AB}| = |\vec{AD}| = |\vec{BC}| = |\vec{CD}| = a$, we compute 
\be  \label{boki}
|\vec{AC}| = 2 a \sin \frac{\psi}{2} \ , \quad |\vec{BD}| = 2 a \sin \frac{\phi}{2} \ .
\ee
Thus all sides of the tetrahedron are expressed in terms of $a, \phi, \psi$. Then
\[
T_1 \nn = \frac{ \vec{AB} \times \vec{AC}  }{ |\vec{AB} \times \vec{AC} |} \ .
\]
Taking into account  \rf{Delta} and \ 
$\nn \perp \vec{AD} \equiv a \Delta_2 \rr$ \ we get
\[
a^2 \Delta_1 \nn\cdot \Delta_2 \rr  = T_1 \nn \cdot \vec{AD} \ .
\]
Similarly (because the triangles $ABC$ and $ADC$ are homothetic) we have 
\[
a^2 \Delta_2 \nn \cdot \Delta_1 \rr = a^2 \Delta_1 \nn\cdot \Delta_2 \rr  = T_1 \nn \cdot \vec{AD} \ .
\]
Finally,
\be \label{HHH}
 \Delta_2 \nn \cdot \Delta_1 \rr = \Delta_1 \nn\cdot \Delta_2 \rr  = \frac{ ( \vec{AB} \times \vec{AC} )  \cdot \vec{AD}  }{ a^2 |\vec{AB} \times \vec{AC} |} = \frac{ \det (\vec{AB} , \vec{AC} , \vec{AD} )}{  a^2 |\vec{AB} \times \vec{AC} |} \ .
\ee
We denote by $H$ the height of the tetrahedron $ABCD$,  perpendicular to the plane $ABC$ (i.e., perpendicular to $\pi (T_1 \rr)$). 
The volume of the tethrahedron $ABCD$ is given by $V_{ABCD} = \frac{1}{3} H P_{ABC}$, and
\be \label{pole}  
P_{ABC} = \frac{1}{2} |\vec{AB} \times \vec{AC} | =  \frac{1}{2} a^2 \sin \psi \ , 
\ee
\[
V_{ABCD} = \frac{1}{6} 
\det (\vec{AB} , \vec{AC} , \vec{AD} ) = \frac{1}{6} \sqrt{\left| \begin{array}{ccc} \vec{AB} \cdot \vec{AB} & \vec{AB} \cdot \vec{AC}  & \vec{AB} \cdot \vec{AD} \\ \vec{AC} \cdot \vec{AB} & \vec{AC} \cdot \vec{AC} & \vec{AC} \cdot \vec{AD} \\ \vec{AD} \cdot \vec{AB} & \vec{AD} \cdot \vec{AC} & \vec{AD} \cdot \vec{AD}  \end{array} \right|} \ .
\]
All entries of the determinant can be expressed by $a, \phi, \psi$ using \rf{boki} and the cosine rule. We get 
\be
\det (\vec{AB} , \vec{AC} , \vec{AD} ) = 4 a^3 \sin \frac{\phi}{2} \sin \frac{\psi}{2} \sqrt{\frac{\cos\phi + \cos\psi}{2}} \ .
\ee
Therefore,
\be \label{HH}
H = \frac{ \det (\vec{AB} , \vec{AC} , \vec{AD} )}{ a^2 \sin \psi} = 
 \frac{4 a \sin \frac{\phi}{2} \sin \frac{\psi}{2}}{\sin\psi} \sqrt{\frac{\cos\phi + \cos\psi}{2}} \ .
\ee
 Then, from \rf{HHH} and \rf{HH}, we have
\be    \label{H}
a^2 \Delta_1 \nn\cdot \Delta_2 \rr  = a^2 \Delta_2 \nn \cdot \Delta_1 \rr = H  \ .
\ee
Now, we express  $K$, given by \rf{Kdisc}, in terms of $a, \phi, \psi$. By virtue of \rf{H} we get
\be  \label{KK}
 a^2 K = - \frac{ H^2 }{a^2 (1 - \cos^2 \phi)} = - \frac{\cos\phi + \cos\psi}{2 \cos^2 \frac{\phi}{2} \cos^2 \frac{\psi}{2}} = \tan^2 \frac{\phi}{2} \tan^2 \frac{\psi}{2} - 1  \ .
\ee
The angle $\theta$, defined in Proposition~\ref{tconst},  can be computed from the triangle $DD'O$, where $O$ is the foot of the height $H$ and $D'$ is the foot of the height of the slant $ABD$. 
The area $P_{ABD}$ is $\frac{1}{2} \sin\phi$, therefore  $|DD'| = \sin\phi$. From Pythagoras' theorem we get (after elementary computations)
\[
 |OD'| = \sqrt{\sin^2 \phi - H^2} = 2 \sin^2 \frac{\phi}{2} \tan \frac{\psi}{2} \ .
\] 
Then $\cos\theta = |OD'|/|DD'|$, which yields
\be \label{teta}
\cos\theta = \tan \frac{\phi}{2} \tan \frac{\psi}{2} \ .
\ee
Comparing \rf{KK} and \rf{teta} we end the proof ($\theta = \const$ by Proposition~\ref{tconst}). 
\end{Proof}

\begin{rem}
The formula \rf{Kdisc} can be considered as a natural discrete analogue of the Gaussian curvature \rf{Kcont}.
\end{rem}

\begin{cor}
The discrete surfaces of discrete Gaussain curvature $K = - 1$ are characterized by the condition $ a = \sin\theta$.  
\end{cor}

The same condition, $d  = \sin\sigma$, appears in the definition of the classical B\"acklund transformation for pseudospherical surfaces \cite{Eis1}. There $d$ is the length of the segment joining a point of a pseudospherical surface and its B\"acklund transform, and $\sigma$ is the angle between the correponding tangent planes. 

\section{Pseudospherical surfaces on time scales}

Corollary~\ref{cordis} shows that the assumptions of Definition~\ref{pseudodis} can be expressed completely in terms of the delta derivatives. First, given an immersion $\rr$ on a time scale, we define the normal vector
\be
\nn := \frac{ D_1 \rr \times D_2 \rr }{\sqrt{ 1 - D_1 \rr \cdot D_2 \rr } } \ .
\ee

\begin{Def} \label{pseudotime} An immersion $\rr: \T_1 \times \T_2 \ni (u,v) \rightarrow \rr (u,v) \in \R^3$ such that for any $u,v \in \T_1 \times \T_2$
\begin{itemize}
\item $\rr$ is completely delta differentiable ,
\item $\nn$ is completely delta differentiable ,
\item $|D_1 \rr| = |D_2 \rr| = 1$ \ ,
\item $D_1 \nn \cdot D_1 \rr =  D_2 \nn \cdot D_2 \rr = 0$ \ ,
\end{itemize}
is called a Chebyshev net on the time scale $\T_1 \times \T_2$ (or a pseudospherical surface on the time scale).
\end{Def}

We conjecture that the Gaussian curvature for pseudospherical surfaces on  time scales is given by the formula analogical to \rf{Kdisc}
\be  \label{Ktime}
 K = - \frac{ (D_1 \nn\cdot D_2 \rr ) (D_2 \nn \cdot D_1 \rr )   }{ 1 - (D_1 \rr \cdot D_2 \rr )^2 }  \ ,
\ee
but the rigorous proof is not availabe yet. 
The formulae \rf{Kcont} and \rf{Kdisc} are particular cases of \rf{Ktime}, when $\T_1 = \T_2 = \R$ and $T_1 = T_2 = a \Z$, respectively.

\section{Conclusions}

In this paper the notion of pseudospherical immersions is extended on the so called time scales,  unifying the continuous and discrete cases in a single framework. In particular, the Gaussian curvature of discrete pseudospherical surfaces is defined in a way admitting a straightforward extension on time scales (Proposition~\ref{Kprop}).
It would be interesting to extend other results of the integrable discrete geometry on time scales. 
This is especially important in the context of the numerical approximation of continuous integrable models.

{\it Acknowledgements.} I am grateful to Iwona \'Swis\l ocka for a cooperation \cite{Sw-mgr} and to Klara Janglajew for turning my attention on references \cite{BG,Hi,Hi2}. My work was partially supported by the KBN grant No.\  1 P03B 017 28.


\begin{thebibliography}{99}

\newcommand{\vbib}{\par \vspace{-2ex} \par \bibitem}

\footnotesize

\vbib{Ampol}
Yu.A.Aminov: 
{\it Contemporary theory of surfaces}, University of Bialystok, Bialystok 2004 [in Polish].


\vbib{BP-pseudo}
A.I.Bobenko, U.Pinkall:
 ``Discrete surfaces with constant negative Gaussian curvature and the 
Hirota equation'', {\it J.\ Diff.\ Geom.} {\bf 43} (1996) 527-611. 

\vbib{BS}
A.I.Bobenko, Yu.B.Suris:
``Discrete differential geometry. Consistency as integrability'', 
{\it preprint:} math.DG/0504358. 

\vbib{BG}
M.Bohner, G.Sh.Guseinov:
``Partial differentiation on time scales'',
{\it Dyn.\ Sys.\ Appl.} 13 (2004) 351-379.

\vbib{Ci-fam}
J.L.Cie\'sli\'nski:
``Discretization of multidimensional 
submanifolds associated with Spin-valued spectral problems'', 
{\it Fund.\ Appl.\ Math.\ } {\bf 12} (1) (2006) 253-262, 
{\it preprint:} nlin.SI/0606010.

\vbib{CDS}
 J.Cie\'sli\'nski,  A.Doliwa, P.M.Santini:
``The Integrable Discrete Analogues of Orthogonal Coordinate Systems are
Multidimensional Circular Lattices'',
{\it Phys.\ Lett.\ } {\bf A 235} (1997) 480-488.



\vbib{Do}
A.Doliwa:
``Integrable multidimensional discrete geometry'',
[in:] {\it Integrable hierarchies and modern physical theories}, 
pp.\ 355-389, 
edited by H.Aratyn and A.S.Sorin, 
Kluwer, Dordrecht 2001.


\vbib{DS-nets}
  A.Doliwa, P.M.Santini:
``Multidimensional quadrilateral lattices are integrable'',
{\it Phys.\ Lett.\ } {\bf A 233} (1997) 365-372.

\vbib{DSM}
A.Doliwa, P.M.Santini, M.Ma\~nas: 
``Transformations of quadrilateral lattices'',
{\it J.\ Math.\ Phys.} {\bf 41} (2000) 944-990.


\vbib{Eis1}
  L.P.Eisenhart:
  {\it A Treatise on the Differential Geometry of Curves and Surfaces},
  Ginn, Boston 1909 (Dover, New York 1960).

\vbib{Hi}
S.Hilger:
``Analysis on measure chains -- a unified approach to continuous and discrete calculus'',
{\it Results Math.} 18 (1990) 18-56.

\vbib{Hi2}
S.Hilger:
``Differential and difference calculus -- unified!'', 
{\it Nonl.\ Anal.\ Theory, Meth.\ Appl.} {\bf 30} (1997) 2683-2694.
             
\vbib{Sau}
  R.Sauer: {\it Differenzengeometrie}, Springer, Berlin 1970
[in German].


\vbib{WKS}
W.K.Schief:
``On the unification of classical and novel integrable surfaces. II. Difference geometry'', 
{\it Proc.\ R.\ Soc.\ London A} {\bf 459} (2003) 373-391.


\vbib{Sym}
   A.Sym: ``Soliton surfaces and their application. Soliton
geometry from spectral problems'', [in:] {\it  Geometric Aspects of the
Einstein Equations and
Integrable Systems} (Lecture Notes in Physics {\bf 239}), 
edited by R.Martini; pp.\ 154-231, Springer, Berlin 1985.


\vbib{Sw-mgr}
I.\'Swis\l ocka: 
``Discretization of pseudospherical surfaces'', 
Master Thesis, University of Bia\l ystok, Faculty of Mathematics and Physics,  Bia\l ystok 2005 [in Polish].

\vbib{Wun}
W.Wunderlich:
``Zur Differencengeometrie der Fl\"achen konstanter negativer Kr\"ummung'',
{\it Sitzungsber.\ Ak.\ Wiss. } {\bf 160} (1951) 39-77 [in German].



\end{thebibliography}
\end{document}